# Semidefinite Relaxation-Based Optimization of Multiple-Input Wireless Power Transfer Systems

Hans-Dieter Lang, *Student Member, IEEE,* and Costas D. Sarris, *Senior Member, IEEE*

*Abstract*—An optimization procedure for multi-transmitter (MISO) wireless power transfer (WPT) systems based on tight semidefinite relaxation (SDR) is presented. This method ensures physical realizability of MISO WPT systems designed via convex optimization — a robust, semi-analytical and intuitive route to optimizing such systems. To that end, the nonconvex constraints requiring that power is *fed into* rather than *drawn from* the system via all transmitter ports are incorporated in a convex semidefinite relaxation, which is efficiently and reliably solvable by dedicated algorithms. A test of the solution then confirms that this modified problem is equivalent (tight relaxation) to the original (nonconvex) one and that the true global optimum has been found. This is a clear advantage over global optimization methods (e.g. genetic algorithms), where convergence to the true global optimum cannot be ensured or tested. Discussions of numerical results yielded by both the closed-form expressions and the refined technique illustrate the importance and practicability of the new method. It, is shown that this technique offers a rigorous optimization framework for a broad range of current and emerging WPT applications.

## I. Introduction

WIRELESS power transfer (WPT) systems have come a long way, since their origins over a century ago [1] and their more recent rediscovery [2]. Among many other advancements and developments, multi-transmitter WPT systems have been investigated, both theoretically and experimentally [3]–[6]. Leveraging the additional design parameters offered by the use of multiple transmitters, these multiple-input single-output (MISO) WPT systems have been shown to outperform their single-input single-output (SISO) counterparts. A key aspect of the improved performance of MISO over SISO WPT systems is the natural resilience of the former to the effects of transmitter-receiver misalignment, which are known to severely compromise the power transfer efficiency (PTE) of SISO WPT systems. On the other hand, the increase in degrees of freedom of such systems also makes it challenging to find their optimum operating mode and estimate their maximum achievable performance. In such cases, global optimization methods such as Genetic Algorithms (GA) are commonly used, see for example [7]; however, there is no guarantee that they converge to the global optimum, despite their significant computational cost.

Recently, the fundamental physical limits of such MISO WPT systems were determined using convex optimization [6]. The maximum achievable PTE as well as all relevant electrical parameters (such as currents, resistive and reactive loading, etc.) for optimal operation were derived in closed form. These closed-form results, which are briefly reviewed in Sec. II, are highly valuable, very intuitive and physically meaningful.

However, it still possible that insufficiently constrained optimization of WPT systems may lead to solutions where power is *drawn from* the system via some transmitter ports, instead of being *fed into* it. Such solutions are mathematically valid optima, whose performance (high PTE) is owed to the "recycling" of power among transmit ports suppressing the total transmit power. Yet, they are undesirable and impractical from a physical point of view.

Hence, the formulation of alternative techniques for the optimization of MISO WPT systems is strongly motivated, where such operating conditions are avoided from the beginning. As will be discussed in detail, adding such transmit power constraints to the original optimization problem is not trivial, as they are nonconvex; thereby breaking convexity of the entire optimization problem. As a result, the asssociated optimization problem is unsolvable in general, as the computational cost grows rapidly when increasing the number of transmitters.

In the following, a MISO WPT optimization strategy based on *tight semidefinite relaxation* is presented, whereby the nonconvex problem is reformulated into a convex form with slightly looser constraints, which can be efficiently and reliably solved by dedicated algorithms. A simple test of the solution then proves that this relaxation is *tight*, i.e. fully equivalent to the original problem and that the true global optimum was found. With this framework, these nonconvex transmit power constraints can be efficiently included in the optimization process, ensuring the practical realizability of the proposed MISO WPT systems. Notably, this new optimization framework is far more general, powerful and versatile than the one previously provided in closed form [6], while retaining its mathematical rigor.

This paper begins with a brief introduction and a review of the closed-form expressions [6]. Then, it is shown in detail how the additional constraints are incorporated into the optimization procedure, how its semidefinite relaxation is derived and how the optimization result can be tested for tightness/equivalence to the original nonconvex problem. Finally, simulation results demonstrate the practicality of the proposed method and the importance of including nonconvex transmit power constraints to obtain practically realizable MISO WPT systems.

Remarks on the notation: Italicized thin letters represent scalars, bold small letters refer to vectors, and bold capital letters are matrices; $\mathbf{v}^T$ stands for the transpose of $\mathbf{v}$, while $\mathbf{v}^H$ stands for its Hermitian (conjugate transpose). Real and imaginary parts of complex quantities are marked by $(\cdot)'$ and $(\cdot)''$, respectively, i.e. $\alpha = \alpha' + j\alpha''$. The complex conjugate is denoted by a bar (overline). The symbols $\succ$ ($\succeq$) and $\prec$ ($\preceq$) are used to denominate positive (semi-) definiteness and negative (semi-) definiteness of matrices, respectively. For vectors, they stand for elementwise positive (nonnegative) and negative (nonpositive) entries. A star $(\cdot)^\star$ marks the optimized arguments that lead to the optimal solution of an optimization problem.

Remarks on the terminology: 'Program' is a synonym for 'optimization problem', commonly used in the context of mathematical optimization [8]. Further, in this paper, 'feasibility' will refer to satisfying constraints within a program; i.e. a feasible solution satisfies all required constraints.

## II. Preliminaries

### A. Power Transfer Efficiency (PTE)

The central figure of merit when optimizing WPT systems is the *power transfer efficiency* (PTE) [3]–[6], [9], [10]: the ratio of the power $P_L$ transferred to the load $R_L$ to the total transmit power provided by the source $P_t$

$$\eta = \frac{P_L}{P_t} = \frac{P_L}{P_l + P_L} = \frac{R_L}{R_l + R_L} \; . \quad (1)$$

$P_l$ is the total power absorbed by the system, due to dissipation and radiation, modeled by the loss resistance $R_l$, as illustrated in Fig. 1.

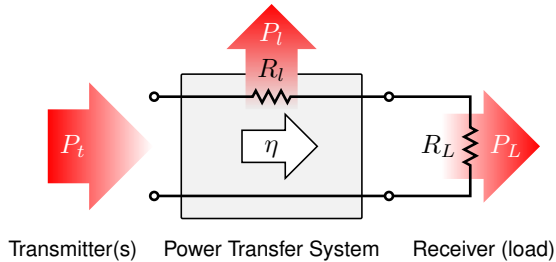

Fig. 1. Illustration of general power transfer systems and the PTE as defined in (1): The ratio of the power transferred to the load, $P_L$, and the total transmitted power, $P_t$. Loss (e.g. due to conduction losses and radiation) is modeled by the resistance $R_l$.

Note that, within the context of this paper, PTE refers only to the electromagnetic power transfer efficiency of the system. Any practical realization of a WPT system will have an end-to-end efficiency which is limited by this PTE; signal generation and impedance matching on the transmitter side as well as matching and rectification will lower the overall performance, but are not considered here. Further, $R_L$ is only practically realizable from a certain resistance level on. However, in this paper the maximum achievable performance in view of the electromagnetic PTE is of central interest, which includes the question of the optimal load resistance.

### B. The MISO WPT System Model

Let the (unloaded) impedance matrix $\mathbf{Z}$ of the MISO WPT system be partitioned according to

$$\mathbf{Z} = \begin{bmatrix} \mathbf{Z}_t & \mathbf{z}_{tr} \\ \mathbf{z}_{tr}^T & z_r \end{bmatrix} \in \mathbb{C}^{N \times N} \; , \quad (2)$$

where the subscripts $t$ and $r$ refer to the transmitter and receiver parts, respectively, and $tr$ stands for the elements coupling the former to the latter.

The diagonal of $\mathbf{Z}_t$ and $z_r$ refer to the loss resistances and self-reactances of the transmitters and receiver, respectively. Typically, $z''_{n,n}, z''_r > 0$ (inductive) when considering systems of magnetically coupled loops or coils. Likewise, the off-diagonal entries of $\mathbf{Z}_t$ and the coupling vectors $\mathbf{z}_{tr}$ contain the mutual inductances $j\omega M_{n,m}$. Generally, each $M_{n,m}$ is complex, due to retardation effects when the electrical distances between the transmitter(s) and/or receiver are not very small.

For physical reasons, impedance matrices have to be positive-real [11]: The real parts of all matrices and submatrices have to be positive-definite, i.e. $\mathbf{Z}', \mathbf{Z}'_t \succ 0$ and $z'_r > 0$. Positive definiteness of a matrix means that the quadratic form of such a matrix is always positive:

$$\mathbf{i}^H \mathbf{Z}' \mathbf{i} > 0 \qquad \forall \mathbf{i} \in \mathbb{C}^N \quad (3)$$

This mathematical property represents the fact that there is always some amount of loss in the system. Positive semidefiniteness ($\succeq 0$) implies nonnegative ($\geq 0$) quadratic forms and for negative (semi-) definiteness the opposite signs and directions apply.

The process of adding reactances to each of the transmitter and receiver nodes as well as a resistive load to the receiver will be referred to as *loading* of the WPT system, where $\hat{\mathbf{Z}}$ is the resulting *loaded impedance matrix*:

$$\hat{\mathbf{Z}} = \mathbf{Z} + j\mathbf{X} + \mathbf{R}_L \; . \quad (4)$$

The real-valued diagonal matrix $\mathbf{X}$ contains the load reactances (positive and negative values referring to inductive and capacitive elements, respectively). The matrix $\mathbf{R}_L$ is zero everywhere but at the last diagonal entry, corresponding to the receiver, where the load resistance $R_L > 0$ is located.

In detail, the voltages and currents $\mathbf{v}, \mathbf{i}$ of the whole (loaded) WPT system are related according to

$$\underbrace{\begin{bmatrix} \mathbf{v}_t \\ v_r = 0 \end{bmatrix}}_{\mathbf{v}} = \Bigg( \underbrace{\begin{bmatrix} \mathbf{Z}_t & \mathbf{z}_{tr} \\ \mathbf{z}_{tr}^T & z_r \end{bmatrix}}_{\mathbf{Z}} + j \underbrace{\begin{bmatrix} \mathbf{X}_t & \mathbf{0} \\ \mathbf{0} & x_r \end{bmatrix}}_{\mathbf{X}} + \underbrace{\begin{bmatrix} \mathbf{0} & \mathbf{0} \\ \mathbf{0} & R_L \end{bmatrix}}_{\mathbf{R}_L} \Bigg)^{\hat{\mathbf{Z}}} \underbrace{\begin{bmatrix} \mathbf{i}_t \\ i_r \end{bmatrix}}_{\mathbf{i}} \quad (5)$$

where, similarly to the impedance matrices, the subscript $t$ and $r$ refer to the transmitter and receiver voltages and currents, respectively. Since the load resistance $R_L$ is actually part of the impedance matrix, KVL states that the corresponding receiver voltage is zero, i.e. $v_r = 0$, as illustrated in Fig. 2.

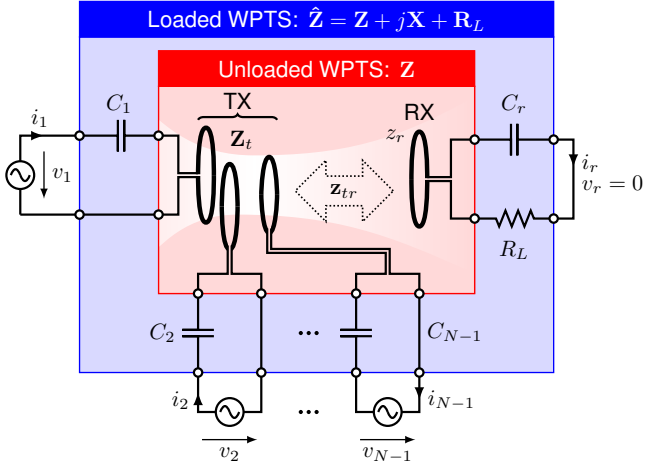

Fig. 2. Loop-based MISO WPT system model: Core structure with unloaded impedance matrix $\mathbf{Z}$ (obtained via simulation or measurements), reactive components, here $x_n = -(\omega C_n)^{-1}$, resistive load $R_L$ and voltage sources $v_n$ added to the loaded impedance matrix $\hat{\mathbf{Z}}$. Note that at the receiver end, since $R_L$ is part of $\hat{\mathbf{Z}}$, the voltage at the corresponding port is zero, $v_r = 0$.

To maximize the performance of such systems, the aim is to find optimal voltages $\mathbf{v}$ and currents $\mathbf{i}$ (real and imaginary parts) as well as loading components $\mathbf{x} = \begin{bmatrix}\mathbf{x}_t^T, x_r\end{bmatrix}^T = \mathrm{diag}(\mathbf{X})$, $R_L$, which maximize the PTE, obtained according to its definition (1) as

$$\eta = \frac{\frac{1}{2}\mathbf{i}^H \mathbf{R}_L \mathbf{i}}{\frac{1}{2}(\mathbf{i}^H \mathbf{v})'} = \frac{\mathbf{i}^H \mathbf{R}_L \mathbf{i}}{\mathbf{i}^H (\mathbf{Z}' + \mathbf{R}_L) \mathbf{i}} \ . \quad (6)$$

As set out above, the PTE can be given in terms of the currents $\mathbf{i}$ alone. The voltages $\mathbf{v}$ follow from (5), but are only physically meaningful as long as $v_r = 0$ is ensured.

### C. Review of the Convex Optimization of MISO WPT Systems

As pointed out previously [6], the biquadratic fraction (6) with both convex numerator and denominator, does not have to be convex itself. However, convexity is the property that ensures that the objective has a single and global optimum that can be reliably found by dedicated algorithms [8].

Note that there are numerical algorithms dedicated to dealing with this type of problem directly, often referring to it as maximization of the *constrained Rayleigh quotient* [12], [13]. However, in the following, a more direct and physically intuitive way of dealing with this problem is taken.

In order to obtain a uniquely defined optimization problem, the receiver current could be assumed to be a known constant $i_r \in \mathbb{R}$, $i_r \neq 0$, (purely real, non-zero). More specifically and without loss in generality, the receiver current is chosen to be $i_r = \sqrt{2/R_L}$, leading to unit transferred power, i.e. $P_L = \frac{1}{2}\mathbf{i}^H \mathbf{R}_L \mathbf{i} = \frac{1}{2}R_L i_r^2 = 1$. Then, maximizing the PTE (6) is equivalent to minimizing the power loss

$$\begin{aligned}P_l &= \frac{1}{2}\mathbf{i}^H \mathbf{Z}' \mathbf{i} \\ &= \frac{1}{2}\left(\mathbf{i}_t'^T \mathbf{Z}_t' \mathbf{i}_t' + \mathbf{i}_t''^T \mathbf{Z}_t' \mathbf{i}_t'' + 2i_r \mathbf{z}_{tr}'^T \mathbf{i}_t' + i_r^2 z_r'\right) \ . \quad (7)\end{aligned}$$

and the resulting PTE is given by

$$\eta_{\max} = \frac{1}{P_{l,\min} + 1} \ . \quad (8)$$

In order to only have to deal with real variables in the following steps, the complex transmitter currents $\mathbf{i}_t$ are separated into real and imaginary parts $\mathbf{i}_t'$ and $\mathbf{i}_t''$, respectively:

$$\mathbf{c}_t = \begin{bmatrix}\mathbf{i}_t' \\ \mathbf{i}_t''\end{bmatrix} \ . \quad (9)$$

To satisfy KVL at the receiver end, i.e. ensuring $v_r = 0$ as in (5) and illustrated in Fig. 2, the following linear constraint (affine [8] in $\mathbf{c}_t$) has to hold:

$$\begin{bmatrix}\mathbf{z}_{tr}'^T & -\mathbf{z}_{tr}''^T \\ \mathbf{z}_{tr}''^T & \mathbf{z}_{tr}'^T\end{bmatrix} \mathbf{c}_t + \sqrt{\frac{2}{R_L}} \begin{bmatrix}\hat{z}_r' \\ \hat{z}_r''\end{bmatrix} = \mathbf{0} \ , \quad (10)$$

where $\hat{z}_r = z_r + R_L + jx_r$ is the loaded receiver self-impedance. Since the receiver reactance $x_r$ is a free parameter, it can always be chosen so as to automatically satisfy the second constraint for the optimal transmitter currents $\mathbf{c}_t^\star$,

$$x_r^\star = -z_r'' - \sqrt{\frac{R_L}{2}} \begin{bmatrix}\mathbf{z}_{tr}''^T & \mathbf{z}_{tr}'^T\end{bmatrix} \mathbf{c}_t^\star \ . \quad (11)$$

Thus, the second row of the constraints (10) can be removed.

Using the objective (7) in terms of the real currents (9) in combination with the remaining constraint of (10), the program to solve for optimal wireless power transfer to a specific load $R_L$ at the receiver can be given in detail as

$$\begin{aligned}P_{l,\min} = \min_{\mathbf{c}_t} \ &\frac{1}{2}\mathbf{c}_t^T \begin{bmatrix}\mathbf{Z}_t' & \mathbf{0} \\ \mathbf{0} & \mathbf{Z}_t'\end{bmatrix} \mathbf{c}_t + \sqrt{\frac{2}{R_L}} \begin{bmatrix}\mathbf{z}_{tr}'^T & \mathbf{0}\end{bmatrix} \mathbf{c}_t + \frac{z_r}{R_L} \\ \text{s.t.} \ &\begin{bmatrix}\mathbf{z}_{tr}'^T & \mathbf{z}_{tr}''^T\end{bmatrix} \mathbf{c}_t = -\sqrt{\frac{2}{R_L}}(z_r' + R_L)\end{aligned} \quad (12)$$

The final optimization problem (12) constitutes a convex *quadratic program* (QP) in terms of the transmitter currents $\mathbf{c}_t$. Such programs can efficiently and reliably be solved to (nearly) arbitrary precision using numerical methods, such as interior-point or active-set algorithms [8]. More importantly however, the fact that (12) is a purely equality-constrained convex QP means that it has an analytic solution, obtainable for example via Lagrangian duality [8]. The steps involved are straightforward and, therefore, omitted here. These closed-form expressions can then be optimized with respect to $R_L$, by setting its derivative to zero. In the following, the so obtained closed-form expressions of [6] are reviewed and presented in a more physically meaningful and intuitive form.

### D. Closed-Form Expressions for Optimal MISO WPT Systems

*1) Fundamental Quantities:* The optimal operating parameters, such as transmitter currents and loading elements, as well as the resulting PTE are found in terms of two fundamental quantities with physical interpretation:

- Minimum-loss output impedance

$$z_o = z_o' + jz_o'' = z_r - \mathbf{z}_{tr}^T \mathbf{Z}_t'^{-1} \mathbf{z}_{tr}' \in \mathbb{C} \quad (13)$$

- Mutual coupling quality factor (also called "mutual $Q$")

$$U = \sqrt{\frac{\mathbf{z}_{tr}^H \mathbf{Z}_t'^{-1} \mathbf{z}_{tr}}{z_o'}} \in \mathbb{R} \quad (14)$$

Since all impedance (sub-) matrices are positive-real, it follows that $z_o', U > 0$.

The real part of the minimum-loss output impedance $z_o$ corresponds to the (minimized) loss due to radiation and dissipation (i.e. $R_l = z_o'$ in (1) and Fig. 1), as illustrated in Fig. 3, whereas the imaginary part $z_o''$ is the reactive component to be canceled out by the receiver reactance $x_r$.

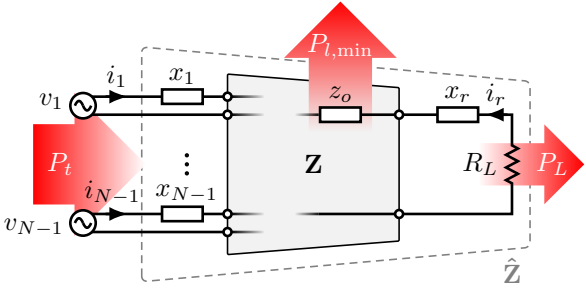

Fig. 3. The optimized MISO WPT system, including the minimal output impedance $z_o$ responsible for both losses during the power transfer as well as output reactance to be compensated.

The presented mutual coupling quality factor is the natural extension of the well-known "mutual $Q$" for SISO systems [9], [14], commonly given as

$$U = k\sqrt{Q_t Q_r} = \frac{\omega M}{\sqrt{R_t R_r}} \quad (15)$$

Note that generally the mutual inductance $M$ and also the coupling coefficient $k$ are complex quantities. In such cases, their absolute values should be used in (15). Evidently, the real part of the output impedance $z_o'$ corresponds to the loss resistance on the receiver side $R_r$, whereas dividing by the loss resistance on the transmitter side $R_t$ is replaced by the inverse of the real part of the transmitter impedance matrix $\mathbf{Z}_t'$. Thus, in terms of quality and coupling factors, the mutual coupling quality factor (14) becomes $U = \sqrt{\mathbf{k}^H \mathbf{Q}_t \mathbf{k} Q_r}$ for MISO systems.

*2) Optimal Currents, Loads and PTE:* The optimal transmitter currents are found as

$$\mathbf{i}_t^\star = -\mathbf{Z}_t'^{-1}\left(\mathbf{z}_{tr}' + \frac{z_o + R_L + jx_r}{z_o' U^2}\overline{\mathbf{z}}_{tr}\right) i_r . \quad (16)$$

The optimal receiver reactance is canceling the reactive component of the output impedance, as mentioned, $x_r^\star = -z_o''$. The transmitter reactances $\mathbf{x}_t = \mathbf{diag}(\mathbf{X}_t)$ are not essential for maximizing the PTE. They can be chosen arbitrarily, as long as the optimal transmitter currents result from applying the voltages to the loaded system.

Further, using the optimal receiver reactance $x_r^\star$, the system operates at the resonant PTE, given by

$$\eta_{\mathrm{res}} = \eta|_{x_r = x_r^\star} = \frac{U^2}{1 + R_L/z_o' + U^2} \cdot \frac{R_L}{R_L + z_o'} . \quad (17)$$

Finally, when in addition to the reactance also using the optimal load resistance

$$R_L^\star = z_o'\sqrt{1 + U^2} \quad (18)$$

the maximum achievable PTE (the physical limit of the WPT system) results:

$$\eta_{\max} = \eta_{\mathrm{res}}|_{R_L = R_L^\star} = \frac{U^2}{(1 + \sqrt{1 + U^2})^2} . \quad (19)$$

Note that $\eta_{\max}$ and $R_L^\star$ are of the same form as for SISO systems [9], [14], where $U$ is given by (15). In other words, SISO WPT systems are special cases of the more general MISO WPT systems, since in the case of a single transmitter, (14) collapses to (15) and real part of (13) is replaced by $R_t$, as previously mentioned; the rest of the expressions remain the same.

### III. PROBLEM STATEMENT

#### A. Port Impedance Matrices (PIMs) $\mathbf{T}_n$

The total transmit power in the denominator of (1) can be separated into the contributions of each port $n$:

$$P_t = P_l + P_L = \frac{1}{2}(\mathbf{i}^H\mathbf{v})' = \frac{1}{2}\mathbf{i}^H\hat{\mathbf{Z}}'\mathbf{i}$$
$$= \sum_n P_{t,n} = \frac{1}{2}(v_n \bar{i}_n)' = \frac{1}{2}\sum_n \mathbf{i}^H \hat{\mathbf{T}}_n \mathbf{i} \quad (20)$$

where $\hat{\mathbf{T}}_n$ are going to be referred to as *port impedance matrices* (PIMs). These PIMs are interesting in many ways:

- They sum up to the total loss resistance matrices:

$$\sum_n \mathbf{T}_n = \mathbf{Z}' \quad \text{and} \quad \sum_n \hat{\mathbf{T}}_n = \hat{\mathbf{Z}}' . \quad (21)$$

 Note that $\mathbf{T}_n = \hat{\mathbf{T}}_n$, for all $n$ except $n = N$, where $\hat{\mathbf{T}}_N = \mathbf{T}_N + \mathbf{R}_L$.
- While $\mathbf{Z}'$ and $\hat{\mathbf{Z}}'$ are purely real-valued, symmetric and positive definite as previously noted, the PIMs are complex-valued, Hermitian, singular and indefinite.
- Each $\mathbf{T}_n$, $\hat{\mathbf{T}}_n$ has exactly one positive, one negative and $N-2$ zero eigenvalues. These eigenvalues and the corresponding eigenvectors can be derived analytically, as laid out in appendix A.

#### B. Nonconvex Transmit Power Constraints

The total transmit power is always $P_t > 0$ (strictly positive), when transferring power $P_L > 0$. However, this is not necessarily true for the transmitter powers at each port $n = 1, ..., N-1$, since all $\mathbf{T}_n$ are indefinite, as previously mentioned.

As will be shown in the results section, in some cases, the optimal currents (16) lead to one or more transmitter powers being negative. This implies that via some transmitter ports, power is *drawn from* the system, rather than *fed into* it. Such negative transmitter powers lead to a reduced net transmit power (20), which has a positive effect on the PTE (1).

While mathematically correct and physically meaningful overall (at least from a circuit theory point of view), such solutions prove difficult to realize in practice, as discussed in more detail in Sec. V. Therefore, it makes sense to try to avoid such solutions altogether, during the optimization process.

In order to accomplish that, the constraints

$$P_{t,n} = \frac{1}{2}\mathbf{i}^H \mathbf{T}_n \mathbf{i} \geq 0 \qquad n = 1, \ldots, N-1 \qquad (22)$$

would have to be added to (12). They ensure nonnegative transmit powers and, thus, power being fed into the system at all transmitter ports, rather than being drawn from it.

Adding the transmitter power constraints (22), the program (12) in terms of complex currents $\mathbf{i}$ becomes:

$$\begin{aligned}
P_{l,\min} = \min_{\mathbf{i}} \ & \frac{1}{2}\mathbf{i}^H \mathbf{Z}' \mathbf{i} \\
\text{s.t.} \ & \frac{1}{2}\mathbf{i}^H \mathbf{T}_n \mathbf{i} \geq 0 \qquad n = 1, \ldots, N-1 \\
& v_r = \begin{bmatrix} \mathbf{z}_{tr}^T, \hat{z}_r \end{bmatrix} \mathbf{i} = 0 \\
& i_r = \begin{bmatrix} \mathbf{0}, 1 \end{bmatrix} \mathbf{i} = \sqrt{\frac{2}{R_L}}
\end{aligned} \qquad (23)$$

Unfortunately, since $\mathbf{T}_n \not\succeq 0$ (indefinite and not negative-semidefinite, as required to make (22) convex), the added constraints are nonconvex, rendering the whole optimization problem nonconvex. Thus, (23) is a *nonconvex QCQP* [8] (as opposed to the convex QCQP mentioned at the end of Sec. II-C). Note that the same would be true for constraints of the form $\frac{1}{2}\mathbf{i}^H \mathbf{T}_n \mathbf{i} \leq P_{t,n,\max}$, i.e. limiting the maximum transmit power at each port.

Nonconvex QCQPs belong to the class of NP-hard problems [15]. This means that there are no known algorithms which can solve the problem for an arbitrary number of transmitters.

## IV. SEMIDEFINITE RELAXATION

In the following, a method is presented which incorporates the nonconvex transmit power constraints while remaining computationally efficiently and reliably solvable.

### A. Nonconvex QCQP

Zero KVL in the imaginary part can always be achieved by appropriately choosing the receiver reactance $x_r$ and the imaginary part of the receiver current is therefore chosen to always remain zero. Thus, the actual problem size is $M = 2N - 1$ and the unknown currents are (in real form)

$$\mathbf{c} = \begin{bmatrix} \mathbf{i}'_t \\ i'_r \\ \mathbf{i}''_t \end{bmatrix} \in \mathbb{R}^M . \qquad (24)$$

The nonconvex QCQP (23) can be formulated in standard QCQP form as follows:

$$\begin{aligned}
P_{l,\min} = \min_{\mathbf{c}} \ & \mathbf{c}^T \mathbf{Q}_0 \mathbf{c} \\
\text{s.t.} \ & \mathbf{c}^T \mathbf{Q}_n \mathbf{c} \geq 0 \qquad n = 1, \ldots, N-1 \\
& \mathbf{A}\mathbf{c} = \mathbf{b}
\end{aligned} \qquad (25)$$

where the objective and the inequality constraints contain the leading principal $M \times M$ submatrices of the real block matrix representations:

$$\mathbf{Q}_n = \frac{1}{2} \begin{bmatrix} \mathbf{T}'_n & -\mathbf{T}''_n \\ \mathbf{T}''_n & \mathbf{T}'_n \end{bmatrix}_{M \times M} \qquad (26a)$$

and $\mathbf{Q}_0 = \sum_n \mathbf{Q}_n$.

Similar to the real part of the impedance matrix and the PIMs, $\mathbf{Q}_0 \succ 0$ while $\mathbf{Q}_n \not\succeq 0$. The equality constraints represent the KVL at the receiver node and the fixed receiver currents (real and imaginary parts). The affine equality constraints $\mathbf{A}\mathbf{c} = \mathbf{b}$ are actually

$$\begin{bmatrix} \mathbf{z}'_{tr} & z_r + R_L & -\mathbf{z}''_{tr} & 0 \\ \mathbf{0} & 1 & \mathbf{0} & 0 \end{bmatrix} \mathbf{c} = \begin{bmatrix} 0 \\ \sqrt{2/R_L} \end{bmatrix} \qquad (27)$$

to ensure KVL at the receiver end ($v_r = 0$) as well as unit transferred power $P_L = 1$, as before.

Using the cyclic property of the trace $\mathbf{c}^T \mathbf{Q}_i \mathbf{c} = \text{tr}(\mathbf{c}^T \mathbf{Q}_i \mathbf{c}) = \text{tr}(\mathbf{Q}_i \mathbf{c}\mathbf{c}^T)$ the objective and all inequality constraints can be written in linear terms of the (quadratic) current matrix $\mathbf{C} = \mathbf{c}\mathbf{c}^T \succeq 0$:

$$\begin{aligned}
P_{l,\min} = \min_{\mathbf{C},\mathbf{c}} \ & \text{tr}(\mathbf{Q}_0 \mathbf{C}) \\
\text{s.t.} \ & \text{tr}(\mathbf{Q}_n \mathbf{C}) \geq 0 \qquad n = 1, \ldots, N-1 \\
& \mathbf{A}\mathbf{c} = \mathbf{b} \\
& \mathbf{C} = \mathbf{c}\mathbf{c}^T
\end{aligned} \qquad (28)$$

Note that the three programs (23), (25) and (28) are mathematically fully equivalent. The non-convexity of the former two has been isolated in the last condition of (28), since the traces are linear operations in terms of the matrices $\mathbf{C}$.

### B. Semidefinite Relaxation in Partially Conic Form

The idea is to exchange the nonconvex constraint with a practically equivalent, yet convex one. Similar to $x = y$ being equivalent to *both* $x \geq y$ *and* $x \leq y$, $\mathbf{C} = \mathbf{c}\mathbf{c}^T$ is equivalent to *both* $\mathbf{C} \succeq \mathbf{c}\mathbf{c}^T$ *and* $\mathbf{C} \preceq \mathbf{c}\mathbf{c}^T$. The *semidefinite relaxation* (SDR) arises from removing the second inequality constraint and using

$$\mathbf{C} - \mathbf{c}\mathbf{c}^T = 0 \quad \xrightarrow{\text{SDR}} \quad \mathbf{C} - \mathbf{c}\mathbf{c}^T \succeq 0 \qquad (29)$$

Implicitly, this also ensures $\mathbf{C} \succeq 0$, since $\mathbf{c}\mathbf{c}^T \succeq 0$.

Via Schur complement, the relaxed constraint (29) is equivalent to

$$\mathbf{C} - \mathbf{c}\mathbf{c}^T \succeq 0 \quad \Leftrightarrow \quad \begin{bmatrix} \mathbf{C} & \mathbf{c} \\ \mathbf{c}^T & 1 \end{bmatrix} \succeq 0 \qquad (30)$$

which is the most common way to express the semidefinite (inequality) constraint. Note that another formulation of the last constraint in (28) would be to require that the rank of the composite matrix in (30) be 1. This implicitly ensures both symmetry and $\mathbf{C} = \mathbf{c}\mathbf{c}^T$, while the equality constraint $\mathbf{A}\mathbf{c} = \mathbf{b}$ ensures a nonzero result.

Thus, the relaxed *semidefinite program* (SDP) of (28), using the semidefinite constraint (30) is

$$P_{l,\min}^{\text{relax}} = \min_{\mathbf{C},\mathbf{c}} \ \text{tr}(\mathbf{Q}_0 \mathbf{C})$$
$$\text{s.t.} \ \text{tr}(\mathbf{Q}_n \mathbf{C}) \geq 0 \quad n = 1, \ldots, N-1$$
$$\mathbf{A}\mathbf{c} = \mathbf{b} \quad (31)$$
$$\begin{bmatrix} \mathbf{C} & \mathbf{c} \\ \mathbf{c}^T & 1 \end{bmatrix} \succeq 0 \ .$$

and represents the dual of the well-known *Shor relaxation* [15] with added affine equality constraints. This problem can be implemented readily in Matlab using CVX [16], [17], a package for specifying and solving convex programs such as this SDP. The actual numerical algorithm used is SDPT3 [18], [19], which for example for a MISO-3 system (problem size $M = 7$) usually converges within ten to thirty iterations, in just a few seconds.

With this relaxed program, a lower bound on the minimum power loss and, thus, by application of (8), an upper bound on the PTE is found: $P_{l,\min}^{\text{relax}} \leq P_{l,\min}$, corresponding to $\eta_{\max}^{\text{relax}} \geq \eta_{\max}$. If equality is achieved, the relaxation is referred to as *tight* [20], referring to the fact that the bounds have zero gap. Evidently, then (31) is a fully equivalent reformulation of (23), (25) and (28) and solving the (convex) SDP is equivalent to solving the original (nonconvex) QCQP. In essence, this also implies that the removed constraint $\mathbf{C} \preceq \mathbf{c}\mathbf{c}^T$ is naturally satisfied by the rest of the program.

Note that, in general, relaxations where the constraints were relaxed (but not the objective) necessarily lead to a remaining feasibility problem, unless they are tight. By removing a constraint, the feasible region for solutions was enlarged and, therefore, the so found solution cannot be feasible for the original problem unless that constraint was loose, i.e. otherwise satisfied to begin with. The remaining feasibility problem is then to find "the next best" solution which is feasible to the original optimization problem but has the smallest gap to the relaxed optimum. Such problems can generally be arbitrarily difficult to solve, as they essentially contain the original nonconvexity and therefore must still be NP-hard.

In this case, a good candidate of a feasible solution is found directly, since the vector $\mathbf{c}^\star$ (optimal solution to $\mathbf{c}$) necessarily satisfies the affine equality constraints (representing the KVL at the receiver end). However, that vector only satisfies the other constraints if the relaxation is tight. Methods to find bounds on these types of inequalities are available, see e.g. [21], but deemed not strict enough to be useful in this case.

As it turns out, while it is intuitive and numerically convenient, the SDP form with a separate vector for the affine constraints (31) is mathematically difficult to deal with, when attempting to prove tightness.

### C. Semidefinite Relaxation in Purely Conic (Quadratic) Form

*1) Quadratic Forms of Affine Constraints:* Instead of defining the receiver current via affine equality constraints, the received power is fixed in its quadratic form directly, i.e. $\frac{1}{2} i_r^2 R_L = 1$ (since $i_r = i_r'$ and $i_r'' = 0$). Let

$$\mathbf{R} = \frac{1}{2} \begin{bmatrix} \mathbf{R}_L & \mathbf{0} \\ \mathbf{0} & \mathbf{R}_L \end{bmatrix}_{M \times M} \quad (32)$$

be the leading principal submatrix of the block form of $\mathbf{R}_L$, similar to (26). Due to the row/column reduction, $\mathbf{R} \succeq 0$ has only one non-zero entry, the $N$th entry on the diagonal, where $R_L/2 \geq 0$ resides. Then, the condition that ensures unit transferred power (and thereby fixes the real part of the receiver current) becomes:

$$\frac{1}{2} i_r^2 R_L = 1 \quad \Leftrightarrow \quad \text{tr}(\mathbf{R}\mathbf{C}) = 1 \ . \quad (33)$$

Further, since the only remaining constraint is homogeneous, it can also be represented in quadratic form, without loss of generality. Let $\mathbf{k} = \begin{bmatrix} \mathbf{z}_{tr}', z_r + R_L, -\mathbf{z}_{tr}'' \end{bmatrix}$ and $\mathbf{K}_0 = \mathbf{k}\mathbf{k}^T$, then

$$\mathbf{k}^T \mathbf{c} = 0 \quad \Leftrightarrow \quad \text{tr}(\mathbf{K}_0 \mathbf{C}) = 0 \quad (34)$$

is responsible to satisfy KVL at the receiver end, i.e. to ensure $v_r = 0$.

For increased numerical accuracy, the following redundant matrix equalities may be added

$$\text{tr}(\mathbf{K}_m \mathbf{C}) = 0 \quad m = 1, \ldots, M \quad (35)$$

where $\mathbf{K}_m = \mathbf{u}_m \mathbf{k}^T + \mathbf{k}\mathbf{u}_m^T$, where $\mathbf{u}_m$, is a unit column vector of zeros everywhere but at the $m$th position. Note that, while $\mathbf{K}_0 \succeq 0$, $\mathbf{K}_m \not\succeq 0$ for all $m = 1, \ldots, M$.

Thus, incorporating the constraints (33) and (34), the semidefinite relaxation of (25) in purely quadratic (conic) form is obtained as:

$$P_{l,\min}^{\text{relax}} = \min_{\mathbf{C}} \ \text{tr}(\mathbf{Q}_0 \mathbf{C})$$
$$\text{s.t.} \ \text{tr}(\mathbf{Q}_n \mathbf{C}) \geq 0 \quad n = 1, \ldots, N-1$$
$$\text{tr}(\mathbf{K}_m \mathbf{C}) = 0 \quad m = 0, \ldots, M \quad (36)$$
$$\text{tr}(\mathbf{R}\mathbf{C}) = 1$$
$$\mathbf{C} \succeq 0$$

Note that the programs (31) and (36) are fully equivalent.

*2) Duality:* The Lagrangian of the semidefinite program (36) is obtained by adding the constraints, weighted by penalty factors, called *dual variables* [8], to the objective:

$$L = \text{tr}\left(\left[\mathbf{Q}_0 - \sum_{n=1}^{N-1} \lambda_n \mathbf{Q}_n - \sum_{m=0}^{M} \nu_m \mathbf{K}_m - \sigma \mathbf{R}\right] \mathbf{C}\right) + \sigma \quad (37a)$$

$$= \text{tr}\left(\left[\mathbf{P}(\boldsymbol{\lambda}, \boldsymbol{\nu}) - \sigma \mathbf{R}\right] \mathbf{C}\right) + \sigma \quad (37b)$$

$$= \text{tr}(\mathbf{Q}\mathbf{C}) + \sigma \quad (37c)$$

where all $\lambda_n \geq 0$ (and therefore the vector $\boldsymbol{\lambda} \succeq 0$), in order to only punish $\text{tr}(\mathbf{Q}_n \mathbf{C}) < 0$. Note that $\mathbf{P} = \mathbf{P}^T$ and thus also $\mathbf{Q} = \mathbf{Q}^T$.

The original problem, when considering duality called the *primal problem*, can be obtained from the Lagrangian by minimizing the objective given that an inner maximization ensures that all constraints are satisfied:

$$P_{l,\min}^{\text{relax}} = \min_{\mathbf{C} \succeq 0} \ \max_{\substack{\boldsymbol{\lambda} \succeq 0 \\ \boldsymbol{\nu}, \sigma}} L \ , \quad (38)$$

If one or more of the constraints are not satisfied, the Lagrangian is unbounded and diverges to $+\infty$. Thus, (38) and (36) are equivalent.

General (weak) duality [8] states that

$$\min_{\mathbf{C}\succeq 0} \max_{\substack{\boldsymbol{\lambda}\succeq 0 \\ \boldsymbol{\nu},\sigma}} L \geq \max_{\substack{\boldsymbol{\lambda}\succeq 0 \\ \boldsymbol{\nu},\sigma}} \min_{\mathbf{C}\succeq 0} L \ . \tag{39}$$

where the right hand side is called the *dual problem*. Since the original problem is convex and optimal feasibility can be presumed for physical reasons, strong duality [8] is expected to hold, leading to equality in (39).

As can easily be seen from (37c), the dual problem is only bounded as long as $\mathbf{Q}$ is positive semidefinite:

$$\min_{\mathbf{C}\succeq 0} \mathrm{tr}(\mathbf{QC}) = \begin{cases} 0 & \mathbf{Q}\succeq 0 \\ -\infty & \text{otherwise} \end{cases} \tag{40}$$

Thus, the dual problem really is

$$P_{l,\min}^{\mathrm{dual}} = P_{l,\min}^{\mathrm{relax}} = \max_{\substack{\boldsymbol{\lambda}\succeq 0 \\ \boldsymbol{\nu},\sigma}} \{\sigma \ : \ \mathbf{Q}\succeq 0\} \tag{41}$$

which is again an SDP, as required (since SDPs are self-dual).

*3) Test of Tightness:* The Karush-Kuhn-Tucker (KKT) conditions [8] of optimality for the two problems (36) and (41) are

$$\mathbf{C}^{\star} \succeq 0 \quad \text{Primal feasibility} \tag{42a}$$
$$\mathrm{tr}(\mathbf{K}_m \mathbf{C}^{\star}) = 0 \quad \text{Primal eq. constraints } 1,\ldots,M \tag{42b}$$
$$\mathrm{tr}(\mathbf{R}\mathbf{C}^{\star}) = 1 \quad \text{Primal equality constraint } M+1 \tag{42c}$$
$$\mathrm{tr}(\mathbf{Q}_n \mathbf{C}^{\star}) \geq 0 \quad \text{Primal ineq. constr. } 1,\ldots,N-1 \tag{42d}$$
$$\boldsymbol{\lambda}^{\star} \succeq 0 \quad \text{Dual feasibility} \tag{42e}$$
$$\mathbf{Q}^{\star} \succeq 0 \quad \text{Dual ineq. (semidef.) constraint} \tag{42f}$$
$$\mathrm{tr}(\mathbf{Q}^{\star}\mathbf{C}^{\star}) = 0 \quad \text{Complementary slackness (CS)} \tag{42g}$$

If the relaxation is tight, $\mathrm{rank}\,\mathbf{C}^{\star} = 1$ or in other words $\mathbf{C}$ has exactly one non-zero eigenvalue and $\mathbf{c}^{\star}$ is the corresponding eigenvector, i.e. $\mathbf{C}^{\star} = \mathbf{c}^{\star}(\mathbf{c}^{\star})^T$. Thus, to test tightness, either the eigenvalues of $\mathbf{C}^{\star}$ can be investigated, or, numerically more efficient, the normalized tightness error, defined as

$$\epsilon = \frac{\left\|\mathbf{C}^{\star} - \mathbf{c}^{\star}(\mathbf{c}^{\star})^T\right\|}{(\mathbf{c}^{\star})^T\mathbf{c}^{\star}} \tag{43}$$

can be used. As long as $\epsilon$ is small (within the bounds of typical numerical approximation errors), this proves that the semidefinite relaxation is tight and, therefore, provides an exact solution $\mathbf{c}^{\star}$ to the nonconvex QCQP (28). Numerical experiments showed that this seems always the case (down to typical normalized tightness errors around $\epsilon \approx 10^{-12}$) and that $\mathbf{C}^{\star}$ has only one nonzero (and indeed positive) eigenvalue and $\mathbf{c}^{\star}$ is its corresponding eigenvector.

### D. Optimizing the Receiver Load $R_L$

The presented method addresses optimization of the transmitter currents and receiver capacitance, but leaves the receiver load resistance $R_L$ unchanged. $R_L$ could be optimized numerically in an outer loop (while optimizing the system for each $R_L$ in the inner loop); as the optimal PTE is a concave function in $R_L$. However, it has been observed that the maximum PTE is usually not sensitive to the load resistance $R_L$, but strongly depends on $x_r \approx x_r^\star$ [22]. Thus, even using the closed-form expression (18), results that are very close to the absolute optimum can be expected, as long as $x_r$ is optimized.

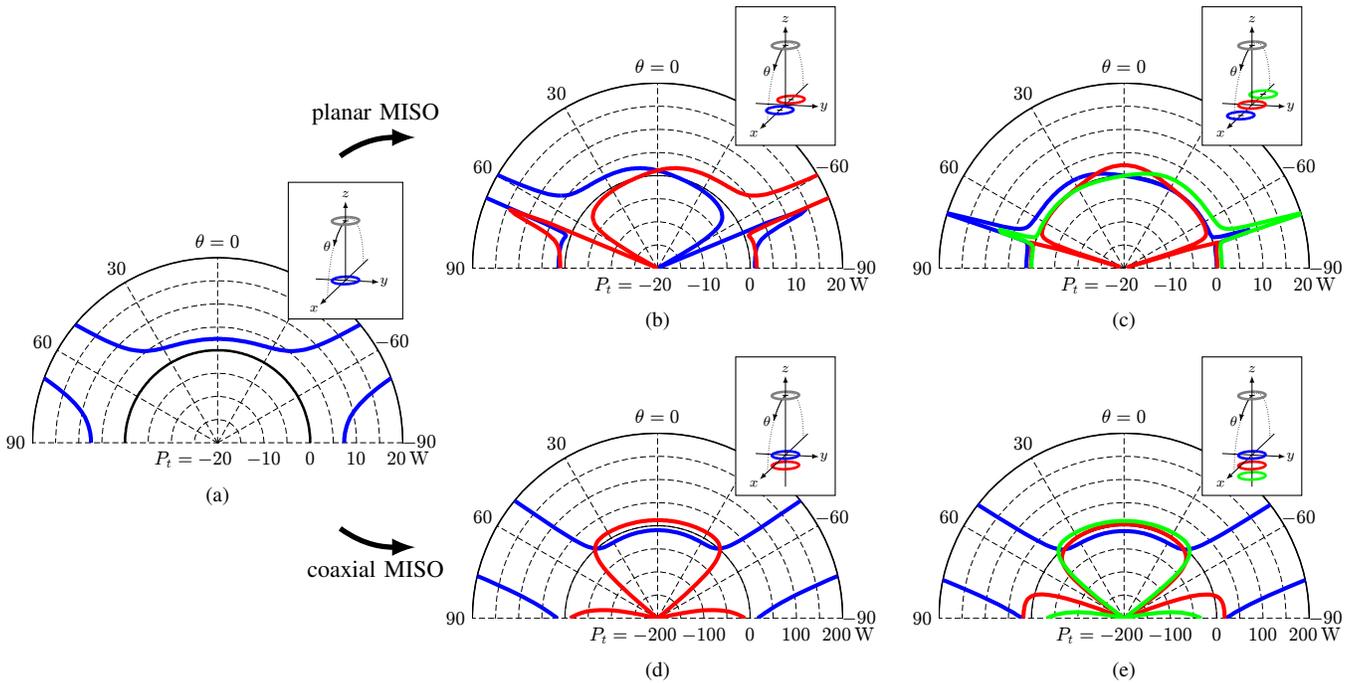

Fig. 7. Comparison of the transmit power patterns of WPT systems optimized via closed-form expressions (13) to (19), [6]: (a) SISO reference, (b) and (c) planar configurations MISO-2p and -3p, respectively, as well as coaxial setups MISO-2c and -3c, (d) and (e), respectively. Negative transmit power mean that power is drawn from, instead of fed into, the system via that particular transmitter port, referenced via the colors. Distance of the receiver $d = \lambda/10$. Note the different scaling of (d) and (e) with respect to all the others.

## V. EXAMPLES

Fig. 4 shows the geometrical setups of the basic SISO and MISO WPT systems under consideration. The frequency of operation is $f = 40\,\text{MHz}$ and the coils are single-turn loops with radii of $r_{\text{loop}} = \lambda/100$, where $\lambda$ is the freespace wavelength at the operating frequency. The conductors are made of copper ($\sigma = 5.8 \times 10^7$ S/m) wire with a thickness (wire radius) of $r_{\text{wire}} = r_{\text{loop}}/10$. The spatial separations of the multiple transmitter loops are $\Delta x = \lambda/50$ and $\Delta z = \lambda/100$, in the $x$- and $z$-directions, respectively. The receiver loop position is in the $xz$-plane ($y = 0$), at a distance $d = 0.05\lambda, \ldots, 0.2\lambda$ from the center of the transmitter (array), as specified for each of the following results, and at the angle $\theta$ off broadside (off the $z$-axis). Note that all the loops are parallel to the $xy$-plane.

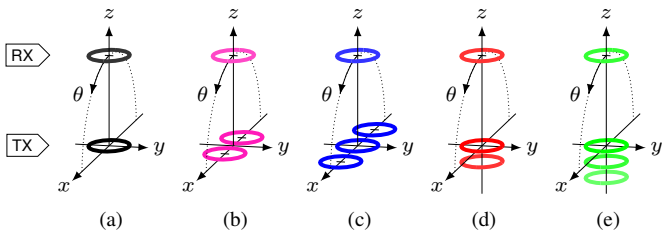

Fig. 4. Illustrations of the SISO WPT reference (a) and the MISO WPT configurations under consideration: Two planar configurations MISO-2p (b) and MISO-3p (c), with transmitters located at the positions $(x_i, y_i, 0)$, and two coaxial configurations MISO-2c (d) and MISO-3c (e) with transmitters at $(0, 0, z_i)$. In all cases a single receiver loop is located at $(d, \theta, \phi = 0)$.

For all the results in this paper, the unloaded impedance matrices were obtained via full-wave simulation using the *Multiradius Bridge-Current* (MBC) method, a computationally efficient and accurate wire-based method of moments (MoM) code with sinusoidal current elements [23], [24].

### A. Closed-Form Optimized MISO WPT Systems

*1) Maximum achievable PTE:* As expected, MISO WPT systems can provide superior performance compared to the SISO reference system: Fig. 5(a) shows typical PTE patterns as a function of the receiver position angle $\theta$ at the distance $d = \lambda/10$ to the transmitter center. Two main observations can be made: First, all the MISO setups outperform the SISO system at every angle $\theta$. However, overall, the increase in PTE is more significant for the planar MISO systems, than for the coaxial ones. This is due to the fact that at angles off broadside ($\theta \neq 0$), the distance to one of the transmitters is always smaller than $d$. Thus, it has to be expected that the transmitter closest to the receiver is favored over the others and that all (or at least most) the power is transferred to the receiver via that particular one. Second, the coaxial setups also provide enhanced performance, predominantly in the broadside direction. In these cases, the transmitters are often being used "in parallel" (as observed when looking at the currents), so as to minimize conduction losses. The performance enhancement appears to become less significant when further increasing the number of transmitters.

*2) Optimal receiver capacitance:* Fig. 5(b) shows the optimal receiver capacitances obtained via the closed-form expressions. It is observed that the capacitance values remain fairly close to about $68.86\,\text{pF}$, over all angles.

As mentioned in the introduction of Sec. IV, previous studies [22] revealed that achieving the maximum possible PTE is generally very sensitive (particularly for such single-turn loop-based MISO setups) to the receiver capacitance values; commonly tolerances much lower than one percent are obtained, to maintain a performance within $1\,\%$ of the optimum PTE. On the other hand, the sensitivity on the load resistance $R_L$ is much smaller; often tolerances of $10\%$ and higher are obtained for the same permissible decrease in performance.

*3) Positive and Negative Transmit Powers:* Fig. 7 shows the underlying working principles of these MISO WPT systems in greater detail, considering the individual transmit powers $P_{t,n}$ as defined in (22): The first graph (a) shows the transmit power of the SISO reference system as a radial function of the receiver position angle $\theta$. As is well known, there are two ranges of angles (around approx. $\pm 60°$ off broadside) at which the transmitter loop couples very weakly to the receiver loop (in fact, infinitely small loops would have points where they do not couple at all). In those areas, the transmitter power spikes to very high values, in order to still transfer unit power ($P_L = 1\,\text{W}$) to the receiver.

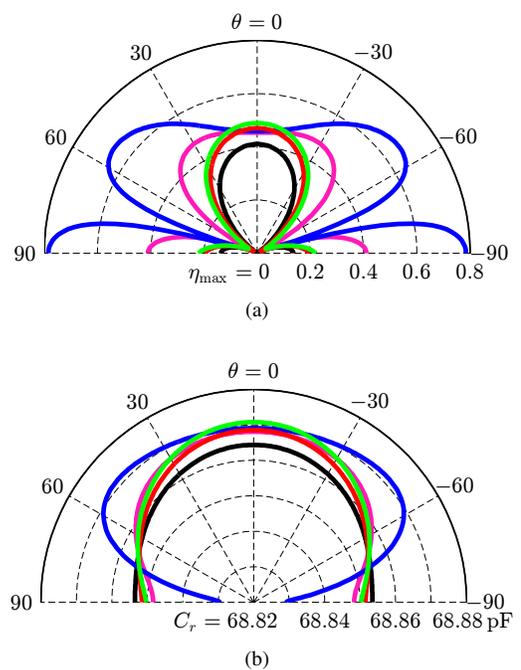

Fig. 5. Comparison of the patterns of the maximum achievable PTE $\eta_{\text{max}}$ (a) and receiver capacitance $C_r^\star$ (b) obtained via closed-form optimization (13) to (19), for the planar and coaxial MISO WPT configurations Fig. 4 (using the same color code) as well as the SISO reference, at the distance $d = \lambda/10$.

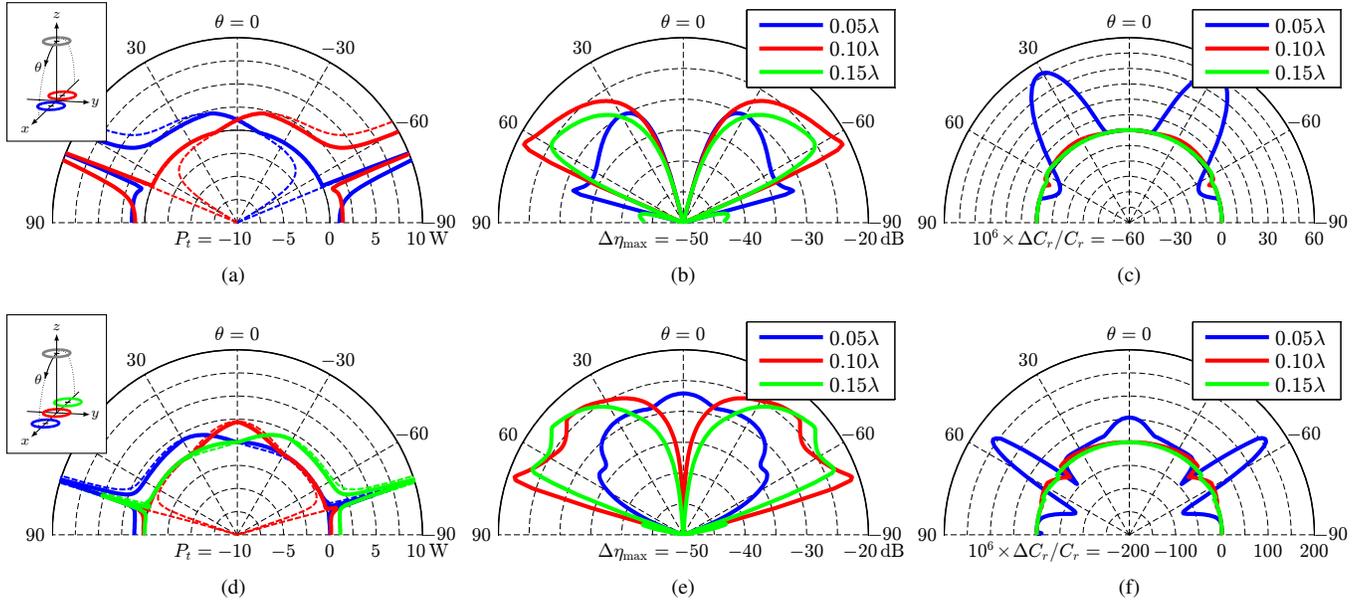

Fig. 9. Transmit powers $P_{t,n}$ of closed-form optima (thin, dashed) and SDR results (thick, solid) at distance $d = 0.1\lambda$ (a, d), PTE reduction with respect to closed-form optimization $\Delta\eta_{\max}$ in dB (b, e), and normalized deviation of receiver capacitance $\Delta C_r/C_r$ (c, f) of the planar multi-transmitter configurations MISO-2p (a-c) and MISO-3p (d-f).

The top right two graphs, in Figs. 7(b) and (c), show the transmit power patterns of the planar MISO setups, whereas (d) and (e) on the bottom right show the coaxial multi-transmitter cases. In the planar cases, most of the power is always transferred via the closest transmitter, as expected. In addition, the remaining transmitters appear to draw away a part of the remaining power in other directions, and feed it back into the system. A similar observation can be made for the coaxial setups: In the regions around the $\pm 60°$ angles, most of the power is inserted into the system by the first (closest) loop, while the others are used to feed some power back into the system. However, at angles close to broadside, most of the power is transferred via the farthest loop, with the remaining loops acting as directors (and also feeding power back into the system). Depending on the distance $d$ between the transmitter center and the receiver, the amplitudes of these negative transmit powers can be quite significant. Note that all loops are modeled with conduction loss and, therefore, extracting power also comes at a cost, unlike in lossless cases; this is already accounted for during the optimization.

### B. MISO WPT Systems with Transmit Power Constraints

Fig. 8 compares one example (MISO-3c, with the color scheme in accordance with the setup shown in Fig. 7(e)), at $\theta = 18°$ off broadside and $d = \lambda/10$ of an optimal solution from the closed-form expressions that involves negative transmit powers (a) and compares it to its nonnegative counterpart (b), to shed light on the difference between the underlying mechanisms. In the first case, power is fed into the system via two transmitter ports, while one transmitter port extracts power. In turn, the net power supplied to the system, given by the sum of all transmitter port powers, is much smaller than the maximum transmit power (but due to the positive-definiteness of the impedance matrix (3) always remains positive). It is implicitly assumed that the power extracted from the system is fully harvested and can be completely "reused", i.e. fed back into the system at 100% efficiency. Evidently, this poses difficulties when trying to realize such a system in practice, as the power cannot be easily "recycled" — at least not at perfect efficiency.

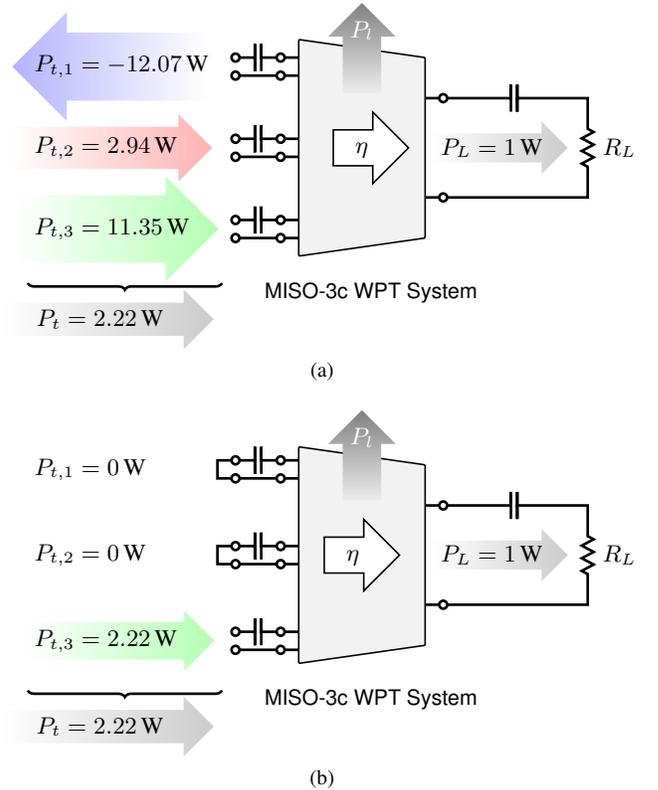

Fig. 8. Examples with and without negative transmit powers: Similar performance (max. PTE of $\eta_{\max} = 1/2.22 \approx 45\%$ in both cases), but case (a) involves harvesting power from one port, while (b) involves only nonnegative transmit powers; the loops connected to the first two ports are only used passively.

The nonnegative transmit power solution in Fig. 8(b) achieves approximately the same performance by only inserting power via the last port and using the optimally tuned transmitter antennas connected to the other two ports passively. This operating mode can be implemented in practice in a straightforward fashion. Additionally, with all transmit powers being positive, they also become smaller in amplitude, which further simplifies the realization of such systems.

Figs. 9 and 11 compare transmit powers, performance degradation and deviation in receiver capacitances in polar form for the two planar and the two coaxial MISO WPT systems, respectively. Fig. 10 shows the corresponding normalized tightness errors, confirming that the relaxation is tight and the solutions solve the original (nonconvex) problem.

The transmit power patterns in the graphs Figs. 9 and 11 (a) and (d) demonstrate that the semidefinite relaxation does indeed lead to solutions where no power is drawn from any transmitter port, as can be seen by comparing the new results (solid thick lines) to the closed-form solutions (thin dashed lines). Note that the optimal constrained transmit power solutions are different from simply setting formerly negative transmit powers to zero (this in fact would both lead to inferior performance as well as break the normalization to unit received power). The performance degradation patterns in the subfigures (b) and (e) reveal that the drop in maximum PTE is minor, in most cases; usually much lower than 1%. The rightmost graphs (c) and (f) show the deviation of the optimal receiver capacitances (see Fig. 5(c)), obtained from the original closed-form expressions (11) and (13). It can be seen that the optimum capacitances remain almost unchanged in all cases, particularly at larger distances.

Lastly, Fig. 10 confirms that the relaxations were tight in all cases, as the normalized tightness errors always remained very small, usually[1] around $10^{-12}$. Gaps in the error plots refer to cases where the closed-form expressions did not lead to any negative transmit powers (no power drawn from the system via any of the transmit ports) and which, therefore, did not have to be optimized using the relaxation method.

## VI. Summary & Conclusion

A closed-form optimization framework [6] for multi-transmitter (MISO) wireless power transfer systems has been reviewed to highlight its simplicity and underlying physical meanings. This powerful and intuitive set of equations is useful for analysis and design of MISO and SISO WPT systems and generalizes previous results for SISO systems [9].

However, as discussed, the obtained absolute optimum operating modes can result in power being drawn via some transmitter ports, rather than fed into the system. Although analytically correct, such operating conditions are difficult to realize in practice as power cannot easily be "recycled".

[1] Larger errors are due to the inherent numerical difficulties of the impedance matrices at larger distances, particularly slight asymmetry.

To ensure practical realizability, it is important to incorporate constraints on the transmit powers into the optimization procedure to avoid such problems from the beginning. However, this is not a trivial task, as such constraints are nonconvex, thereby rendering the entire optimization problem nonconvex.

This paper presents a step-by-step derivation of a convex semidefinite relaxation of the nonconvex problem, where by loosening some constraints a convex formulation is obtained, which can efficiently and reliably be solved numerically by dedicated algorithms. Furthermore, a simple test is given, which reveals whether the solved problem is equivalent to the original nonconvex problem (i.e. whether the relaxation is tight) and that resulting optimum is indeed the true global optimum looked for. This is a clear advantage over global optimization methods, such as genetic algorithms (GA), where convergence to the true global optimum cannot be guaranteed or tested.

Simulation results of some basic multi-transmitter WPT systems show that the new method is very practical and indeed able to avoid any realizability issues due to negative transmit powers. Furthermore, it is shown that the resulting performance gap between the closed-form (absolute) optimum and the new, practically realizable solution is usually very small; typically below 1% of the predicted PTE. Lastly, plots of the tightness test results reveal that in all cases under consideration the relaxation is indeed tight, i.e. that the solved problem is equivalent to the actual nonconvex problem and that the obtained optimum PTE is the global optimum PTE satisfying all constraints.

This new optimization procedure is powerful, versatile and retains the rigor of convex optimization, while ensuring the practical realizability of the optima it produces.

## Appendix A
## Analytical Formulation of the Positive and Negative PIM Parts by Eigenvectors

Let the impedance matrix of the unloaded WPT system be

$$\mathbf{Z} = \mathbf{R} + j\omega\mathbf{L} + j\omega\mathbf{M} \, , \tag{44}$$

where $\mathbf{R}$ and $\mathbf{L}$ are real-valued diagonal matrices containing the losses and self-reactances (usually inductances $L_n > 0$, for loop-based magnetically coupled systems) of the transmitters and the receiver. $\mathbf{M} = \mathbf{M}^T$ is a symmetric (due to reciprocity of the passive system) hollow matrix containing the mutual impedances $j\omega M_{n,m}$. Generally, each $M_{n,m}$ is complex-valued, due to retardation effects when the electrical distance between the transmitter(s) and receiver are not very small.

Moreover, let

$$\{\lambda_{n,m}, \mathbf{v}_{n,m}\} = \mathrm{eig}(\mathbf{T}_n) \qquad n = 1, ..., N \tag{45}$$

denote the $m$th eigenvalue $\lambda_{n,m}$ and corresponding eigenvector $\mathbf{v}_{n,m}$ of the $N \times N$ port impedance matrix (PIM) $\mathbf{T}_n$. In

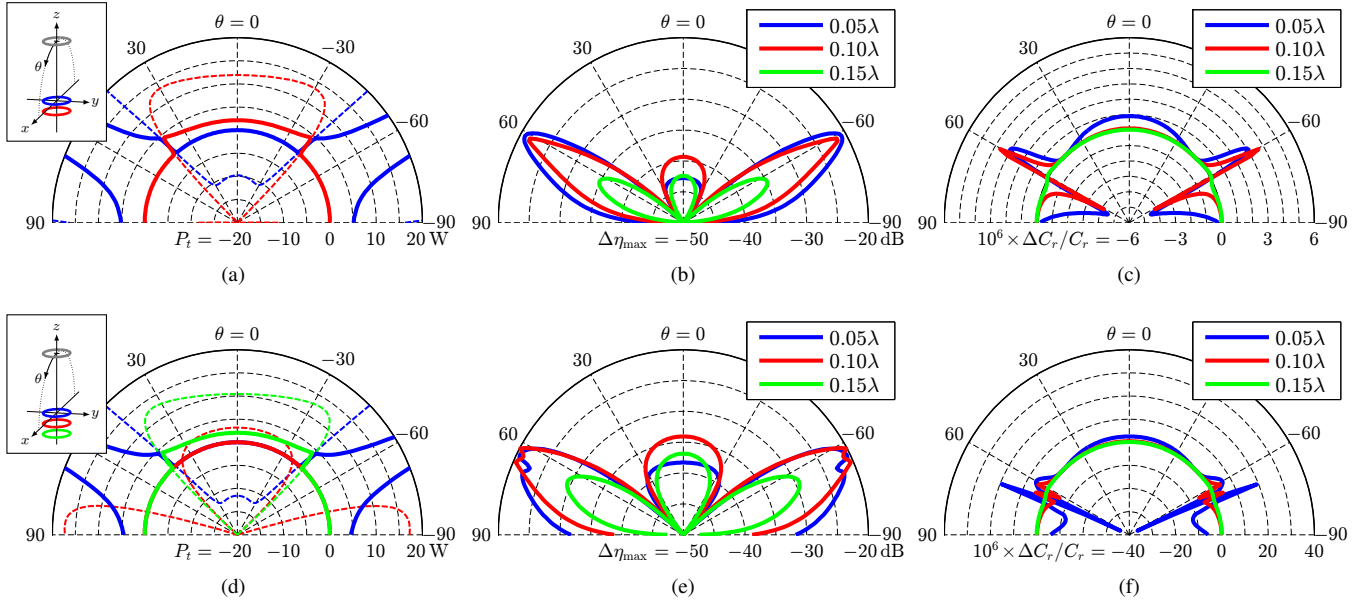

Fig. 11. Transmit powers $P_{t,n}$ of closed-form optima (thin, dashed) and SDR results (thick, solid) at distance $d = 0.1\lambda$ (a, d), PTE reduction to closed-form optimization $\Delta\eta_{\max}$ in dB (b, e), and normalized deviation of receiver capacitance $\Delta C_r/C_r$ (c, f) of the coaxial multi-transmitter configurations MISO-2c (a-c) and MISO-3c (d-f).

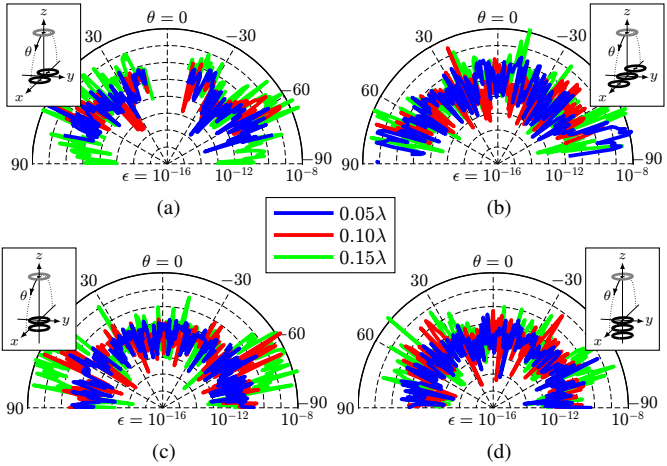

Fig. 10. Normalized tightness errors of the optimal SDR solutions of the four test cases, corresponding to Figs. 9 and 11. Apart from a few exceptions, the erros are usually in the range of $10^{-12}$, confirming that the solution vectors represent the full semidefinite solution well and the relaxation is tight.

the order "negative, positive, and zero", the eigenvalues can be given as

$$\lambda_{n,m} = \begin{cases} -\lambda_n^- = \frac{1}{2}\left(R_n - \sqrt{S_n^2 + R_n^2}\right) < 0 & m = 1 \\ \lambda_n^+ = \frac{1}{2}\left(R_n + \sqrt{S_n^2 + R_n^2}\right) > 0 & m = 2 \\ 0 & m = 3,...,N \end{cases} \quad (46)$$

where the shorthand $S_n^2 = \omega^2 \sum_m |M_{n,m}|^2 > 0$ has been used. Note that, with these definitions, $\lambda_n^+, \lambda_n^- > 0$, for all $n$.

Finally, let the eigenvectors corresponding to the non-zero eigenvalues be denoted by $\mathbf{v}_{n,m=1,2} = \mathbf{v}_n^{\mp}$ (omitting the index $m$, as it is clear from (46) that $m = 1, 2$ correspond to the superscripts $-, +$, respectively) and a superscript zero point to eigenvectors corresponding to zero eigenvalues, i.e. $\mathbf{v}_{n,m>2} = \mathbf{v}_{n,m}^0$ (the index $m$ starts at 3, for these eigenvectors). Hence, the quadratic forms with respect to the eigenvalues and eigenvectors are given by

$$(\mathbf{v}_n^{\pm})^H \mathbf{T}_n \mathbf{v}_n^{\pm} = \pm\lambda_n^{\pm} \cdot (\mathbf{v}_n^{\pm})^H \mathbf{v}_n^{\pm} \quad (47a)$$
$$\mathbf{v}_{n,m}^H \mathbf{T}_n \mathbf{v}_{n,m} = 0 \quad \forall m \neq 1, 2. \quad (47b)$$

The PIMs can be separated into their positive and negative (semidefinite) parts

$$\mathbf{T}_n = \mathbf{T}_n^+ - \mathbf{T}_n^-, \quad (48)$$

where both $\mathbf{T}_n^+, \mathbf{T}_n^- \succeq 0$. Further, each part is simply obtained from its eigenvalues and eigenvectors

$$\mathbf{T}_n^{\pm} = \lambda_n^{\pm} \frac{\mathbf{v}_n^{\pm}(\mathbf{v}_n^{\pm})^H}{(\mathbf{v}_n^{\pm})^H \mathbf{v}_n^{\pm}} \quad (49)$$

with the denominator being the outer product and the numerator the inner product of the respective eigenvectors.

These eigenvectors $\mathbf{v}_n^{\pm}$, can be obtained analytically as well:

$$\mathbf{v}_n^{\pm} = \frac{1}{\tilde{M}_{n,N}^{\pm}} \begin{bmatrix} \tilde{M}_{n,1}^{\pm} \\ \vdots \\ \tilde{M}_{n,N}^{\pm} \end{bmatrix} \quad (50)$$

where $\tilde{M}_{n,m}^{\pm}$ are the entries of the matrix $\tilde{\mathbf{M}}$ which is identical to the mutual impedance matrix $\mathbf{M}$, as given in (44), with the exception of the diagonal:

$$\tilde{\mathbf{M}}^{\pm} = \mathbf{M} \pm 2\mathbf{Diag}(\boldsymbol{\lambda}^{\pm}). \quad (51)$$

Obtaining the positive and negative parts of the PIMs analytically, directly from the impedance matrix entries, adds both computational efficiency as well as numerical precision as compared to using a numerical eigenvalue decomposition.